\documentclass[12pt,intlim,righttag]{article}
\usepackage{amsmath,amsthm}
\usepackage{amssymb}
\usepackage{amsfonts}

\newcommand{\la}{\langle}
\newcommand{\ra}{\rangle}

\newcommand{\sg}{\sigma}

\newcommand{\vp}{\varphi}

\newcommand{\lbl}{\label}

\newcommand\beq{\begin{equation}}
\newcommand\eeq{\end{equation}}

\newcommand{\bea}{\begin{eqnarray}}
\newcommand{\eea}{\end{eqnarray}}
\newcommand{\beaa}{\begin{eqnarray*}}
\newcommand{\eeaa}{\end{eqnarray*}}

\theoremstyle{Theorem}

\theoremstyle{corollary}

\theoremstyle{remark}

\theoremstyle{definition}

\def\a{\alpha}

\def\f{\varphi}

\def\f{\varphi}

\begin{document}
\title{  Martingale Transformations of Brownian Motion with Application  to Functional Equations
 }

\author{M. Mania$^{1)}$ and R. Tevzadze$^{2)}$}

\date{~}
\maketitle

\begin{center}
$^{1)}$ A. Razmadze Mathematical Institute of Tbilisi State University  and
Georgian-American University, Tbilisi, Georgia,
\newline(e-mail: misha.mania@gmail.com)
\\
$^{2)}$ Georgian-American University and
Institute of Cybernetics of Georgian Technical Univercity,  Tbilisi, Georgia
\newline(e-mail: rtevzadze@gmail.com)
\end{center}

\begin{abstract}
{\bf Abstract.}

We describe the classes of functions $f=(f(x), x\in R)$, for which processes $f(W_t)-Ef(W_t)$ and
$f(W_t)/Ef(W_t)$ are martingales. We apply these results to give a martingale characterization of general solutions of
the quadratic and  the D'Alembert  functional equations. We study also the time-dependent martingale transformations
of a Brownian Motion.
\end{abstract}

\noindent {\it 2010 Mathematics Subject Classification. 60G44, 60J65, 97I70}

\noindent {\it Keywords}:  Martingales,  Brownian Motion,  Functional Equations.

\

\section {Introduction.}
It is well known that if for a function $f=(f(x), x\in R)$ the transformed  process $(f(W_t), t\ge0)$  of a Brownian Motion $W$
is a right-continuous martingale, then $f$ is a linear function (see, Theorem 5.5 from \cite{CI}). It is also known that the time-dependent function $f=(f(t,x), t\ge0,x\in R)$ is a linear function of $x$
if and only if the transformed process $(f(t, M_t), t\ge0)$ is a martingale for any martingale $M$  and to this end  to require the martingale property of $f(t,\sigma W_t)$  for two different
$\sigma_1\neq\sigma_2\neq0$ is sufficient (see Corollary 1 of Theorem 2 from \cite{MT}).

In this paper we give  simple generalizations of these results.
We describe the classes of functions $f$ for which the processes $f(W_t)-Ef(W_t)$ and
$f(W_t)/Ef(W_t)$ (for $f(x)>0$) are martingales. We prove that the process
$$
f(W_t)-Ef(W_t)\;\;\;\;\;(\text{resp.}\;\;\;f(W_t)/Ef(W_t)), \;\;\;\;t\ge0
$$
is a  right-continuous  martingale if and only if the function $f(x)$ is of the form
$$
ax^2+bx+c\;\;\;(\text{resp.}\;\;\; a e^{\lambda x}+b e^{-\lambda x}),
$$
for some constants $a, b, c$ and $\lambda$.  Besides, we show that if $f(W_t)-Ef(W_t)$ (resp. $f(W_t)/Ef(W_t))$  is only a martingale (without assuming the regularity of paths), then
$f(x)$ is equal to some square trinomial (resp. to the function $a e^{\lambda x}+b e^{-\lambda x}$ for some $a,b,\lambda$) almost everywhere with respect to the Lebesgue measure.

Our main motivation to consider such martingale transformations of a Brownian Motion was their relations with functional equations. We show that if the function
$f=(f(x), x\in R)$ is a measurable solution of the quadratic functional equation
\begin{equation}\label{i1}
f\left(x+y\right)+f\left(x-y\right)=2f\left(x\right)+2f\left(y\right)\;\;\;\text{for all}\;\;\;x,y \in R,
\end{equation}
then the difference $f(W_t) - Ef(W_t)$ is a martingale and if $f$ is a strictly positive
 solution of the D'Alembert functional equation
\begin{equation}\label{i2}
f\left(x+y\right)+f\left(x-y\right)=2f\left(x\right) f\left(y\right)\;\;\;\text{for all}\;\;\;x,y \in R,
\end{equation}
then a martingale will be the process  $f(W_t)/Ef(W_t)$. The above-mentioned  descriptions of martingale functions enable us to
give equivalent  characterization  of the general measurable solution of equations (\ref{i1}) and (\ref{i2}) in martingale terms.

We consider also time-dependent functions $(f(t,x), t\ge0,x\in R)$ for which the transformed  processes
\begin{equation}\label{i3}
{f(t, \sigma W_t)}-{Ef(t,\sigma W_t)}\;\;\;\text{ and}\;\;\;{f(t, \sigma W_t)}/{Ef(t,\sigma W_t)}
\end{equation}
are martingales, where $\sigma$ is a constant. To obtain simple structural properties for such functions, as for the case of functions $f=(f(x), x\in R)$, one needs some type
of growth conditions on the function $f$, or one should  require the martingale property for transformed processes (\ref{i3})  at least for two different $\sigma\neq 0$. Corresponding assertions (Theorems 5-7)
are given in section 4.

\

\section{Martingale transformations of a Brownian Motion}

Let $W=(W_t, t\ge 0)$ be a standard Brownian Motion defined on a  complete probability space  $(\Omega, {\cal F}, P)$ with filtration
$F=({{\cal F}}_t,t\ge0)$  generated by the Brownian Motion $W$.

{\bf Theorem 1}. Let $f=(f(x), x\in R)$ be a measurable strictly positive  function . Then

a) $f(W_t)$ is integrable for every $t\ge0$  and the  process
$$
N_t=\frac{f(W_t)}{Ef(W_t)},\;\;\; t\ge0,
$$
is a   right-continuous ($P$-a.s.) martingale if and only if the function $f$ is of the form
\begin{equation}\label{form}
f(x)=a e^{\lambda x}+b e^{-\lambda x},
\end{equation}
for some constants $a\ge0,b\ge0, ab\neq 0$ and $\lambda\in R$.

b)  $f(W_t)$ is integrable for every $t\ge0$ and the  process
$N_t$ is a martingale  if and only if the function  $f(x)$  coincides with the function $a e^{\lambda x}+b e^{-\lambda x}$
(for some constants  $a\ge0,b\ge0, ab\neq 0$ and $\lambda\in R$) almost everywhere with respect to the Lebesgue measure.
\begin{proof}
a) Let the process $N_t$ is a   right-continuous ($P$-a.s.) martingale and let
$$
g(t)\equiv Ef(W_t)=\int_Rf(y)\frac{1}{\sqrt{2\pi t}}e^{-\frac{y^2}{2t}}dy.
$$
Since $E|f(W_t)|<\infty$ for all $t\ge 0$, the function $g(t)$ will be continuous for  any $t>0$.
Since $N_t$ is right-continuous  ($P$-a.s.)
and $g(t)$ is continuous, the process $f(W_t)$ will be also right-continuous $P$-a.s..
This implies that  the function $f(x)$ is  continuous  (see Theorem 5.5 from \cite{CI} or Lemma A1 from \cite{MT2}).

Let
$$
F(t,x)=\frac{f(x)}{g(t)}, \;\;\;\; t\ge 0, x\in R.
$$
Since $F(t, W_t)$ is a martingale, we have that
\begin{equation}\label{n1}
F(t,W_t)=\frac{f(W_t)}{g(t)}=\frac{1}{g(T)}E(f(W_T)/{\cal F}_t)
\end{equation}
$P-a.s.$ for all $t\le T$. Let
$$
u(t,x)=E(f(W_T)/ W_t=x).
$$
Since $f$ is positive, $u(t,x)$ will be of the class $C^{1.2}$ on $(0,T)\times R$ and  satisfies the ''backward'' heat equation (see, e.g. \cite{KAR} page 257)
\begin{equation}\label{kol}
\frac{\partial u}{\partial t}+\frac{1}{2}\frac{\partial^2u}{\partial x^2}=0,\;\;\;\;0<t<T, x\in R.
\end{equation}
By the Markov property of $W$
$$
u(t,W_t)=E(f(W_T)/{\cal F}_t)
$$
and  from (\ref{n1}) we have that
$$
f(W_t)=\frac{g(t)}{g(T)}u(t,W_t)\;\;\;\;a.s.
$$
Therefore,
$$
\int_R{\big|f\left(x\right)-\frac{g(t)}{g(T)}u\left(t, x\right)\big|\frac{1}{\sqrt{2\pi t}}e^{-\frac{x^{2}}{2t}}dx}=0
$$
which implies that for any $0<t\le T$
\begin{equation}\label{con}
f\left(x\right)=\frac{g(t)}{g(T)}u\left(t, x\right) \quad
\end{equation}
almost everywhere with respect to the Lebesgue measure. Since  $f$ and $u$ are continuous, we obtain that for any $0<t<T$
equality (\ref{con}) is satisfied for all $x\in R$.

Since $g(t)>0$ for all $t$, this implies that  $g(t)$ is differentiable, $f(x)$ is two-times differentiable and
\begin{equation}\label{ugf}
u\left(t, x\right)=\frac{g(T)}{g(t)}f\left(x\right) \;\;\;\text{for all}\;\;\; x\in R
\end{equation}
for any $0<t<T$.

Therefore, it follows from (\ref{kol}) and (\ref{ugf}) that
$$
\frac{1}{2}\frac{g(T)}{g(t)}f''(x)- \frac{g(T)g'(t)}{g^2(t)}f(x)=0,
$$
which implies that
\begin{equation}\label{fg}
\frac{f''(x)}{f(x)}=2\frac{g'(t)}{g(t)}.
\end{equation}
Since the left-hand side of (\ref{fg}) does not depend on $t$ and the right-hand side on $x$,  both parts of (\ref{fg}) are equal to a constant, which
we denote by $c$. If $c<0$, then the general solution of equation $f''(x)=cf(x)$ leads to $f(x)$ which, with necessity, changes its sign,
hence $c=\lambda^2$ for some $\lambda\in R$. Therefore, we obtain that
$$
f''(x)=\lambda^2 f(x)\;\;\;\text{and}\;\;\;g'(t)=\frac{\lambda^2}{2}g(t).
$$
for some constant $\lambda\in R$. Thus,
$$
f(x)=a e^{\lambda x}+b e^{-\lambda x},\;\;\;g(t)=Ef(W_t)=(a+b)e^{\frac{\lambda^2}{2}t}.
$$
Since the function $f$ should be strictly positive, we shall have that  $a\ge0,b\ge0, ab\neq 0$

Now let us assume that the function $f$ is of the form (\ref{form}). Then $g(t)=E(a e^{\lambda W_t}+ b e^{-\lambda W_t})=
(a+b)e^{\frac{\lambda^2}{2}t}$ and the process
$$
\frac{f(W_t)}{Ef(W_t)}= \frac{a}{a+b}e^{\lambda W_t-\frac{\lambda^2}{2}t}+ \frac{b}{a+b}e^{-\lambda W_t-\frac{\lambda^2}{2}t}
$$
is  a martingale, as a linear combination of two exponential martingales.

b) Let $N_t$ be a martingale  and let
$$
\tilde f(x)=\frac{g(t_0)}{g(T)}u(t_0,x),
$$
for some $t_0>0$. It follows from (\ref{con}) that
\begin{equation}\label{leb}
\lambda (x: f(x)\neq \tilde f(x))=0,
\end{equation}
where $\lambda$ is the Lebesgue measure and  by definition of $u(t,x)$ the function $\tilde f(x)$ is continuous (moreover, it is two-times differentiable).
It follows from (\ref{leb}) that $P(f(W_t)=\tilde f(W_t))=1$ for any $t\ge0$ and since $Ef(W_t)=E\tilde f(W_t)$,
we obtain that for any $t\ge0$
$$
P\big( \frac{f(W_t)}{Ef(W_t)}=\frac{\tilde f(W_t)}{E\tilde f(W_t)}\big)=1.
$$
This implies that the process $\tilde f(W_t)/E\tilde f(W_t)$
is a continuous martingale and it follows from part a) of this theorem that $\tilde f(x)$ is of the form (\ref{form}).
Therefore,
$f(x)$ coincides with  the function $a e^{\lambda x}+b e^{-\lambda x}$ almost everywhere with respect to the Lebesgue measure.

The converse is proved similarly to the part a) of this theorem.
\end{proof}

{\bf Theorem 2}. Let $f=(f(x), x\in R)$ be a  measurable  function , such that  $f(W_t)$ is integrable for every $t\ge0$. Then

a) the  process
$$
M_t={f(W_t)}-{Ef(W_t)},\;\;\; t\ge0,
$$
is a  right-continuous ($P$- a.s.) martingale if and only if the function $f$ is of the form
\begin{equation}\label{form23}
f(x)=a x^2+bx+c \;\;\;\;\text{for some constants}\;\;\;\a, b\;\;\;\text{and}\;\;\;c\in R,
\end{equation}
b)  the process $M_t$ is a martingale if and only if  $f(x)$ coincides with the function $a x^2+bx+c$ (for some constants
$a, b, c\in R$) almost everywhere with respect to the Lebesgue measure.

\begin{proof}
a) Let  the  process
$M_t$ be  a  right-continuous ($P$- a.s.) martingale and let
$g(t)\equiv Ef(W_t)=\int_Rf(y)\frac{1}{\sqrt{2\pi t}}e^{-\frac{y^2}{2t}}dy. $ Using the same arguments as in the proof of Theorem 1,
 the process $f(W_t)$ will be also right-continuous $P$-a.s., which
implies that  the function $f(x)$ is  continuous.

 Let
$$
F(t,x)={f(x)}- {g(t)}, \;\;\;\; t\ge 0, x\in R.
$$
Since $F(t, W_t)$ is a martingale, we have that
\begin{equation}\label{m1}
F(t,W_t)={f(W_t)}-{g(t)}=E(f(W_T)/{\cal F}_t)-g(T)
\end{equation}
$P-a.s.$ for all $t\le T$. Let
$$
u(t,x)=E(f(W_T)/ W_t=x).
$$
Since
\begin{equation}\label{int5}
 E|f(W_t)|=\int_R|f(y)|\frac{1}{\sqrt{2\pi t}}e^{-\frac{y^2}{2t}}dy<\infty \;\;\;\text{for all}\;\;\;t\ge0,
 \end{equation}
  $u(t,x)$ will be of the class $C^{1.2}$ on $(0,T)\times R$ and  satisfies the  ''backward'' heat  equation
 (see, e.g. \cite{KAR})
\begin{equation}\label{mkol}
\frac{\partial u}{\partial t}+\frac{1}{2}\frac{\partial^2u}{\partial x^2}=0,\;\;\;\;0<t<T, x\in R.
\end{equation}
Note that, in Theorem 3.6  of Karatzas and Shreve, using the same arguments, the positivity assumption
on the function $f$
by  integrability condition  (\ref{int5}) can be replaced, which is also sufficient   to guarantee that  (\ref{mkol}) is satisfied.

Similarly to the proof of Theorem 1 one can show that
 for any $0<t\le T$
\begin{equation}\label{mcon}
f\left(x\right)=u\left(t, x\right) -g(T)+g(t)\quad
\end{equation}
almost everywhere with respect to the Lebesgue measure. By continuity of  $f$ and $u$  we obtain that for any $0<t\le T$
the equality (\ref{mcon}) is satisfied for all $x\in R$.

This implies that  $g(t)$ is differentiable, $f(x)$ is two-times differentiable and
 it follows from (\ref{mkol}) and (\ref{mcon}) that
\begin{equation}\label{mf2}
\frac{1}{2}f''(x)=g'(t).
\end{equation}
Since the left-hand side of (\ref{mf2}) does not depend on $t$ and the right-hand side on $x$,  both parts of (\ref{mf2}) are equal to a constant.
Therefore, we obtain that
\begin{equation}\label{f3}
f''(x)=2a \;\;\;\;\text{and}\;\;\; g'(t)=a\;\;\;\;\text{for some}\;\;\;\;a\in R.
\end{equation}
The general solutions of these equations are
\begin{equation}\label{f4}
f(x)=ax^2+bx+c\;\;\;\;\text{and}\;\;\;g(t)=at +c
\end{equation}
respectively, for some $a, b, c\in R$.

Conversely, if  the function $f$ is of the form (\ref{form23})
then
$$
f(W_t)= aW_t^2+bW_t+c,\;\;\;\;Ef(W_t)=at+c
$$
and the process $f(W_t)-Ef(W_t)=a(W_t^2-t)$ is a martingale.

b) Let
$$
\tilde f(x)= u(t_0,x)+g(t_0)-g(T).
$$
for some $t_0>0$. It follows from (\ref{mcon}) that
\begin{equation}\label{leb2}
\lambda (x: f(x)\neq \tilde f(x))=0,
\end{equation}
where $\lambda$ is the Lebesgue measure and  by definition of $u(t,x)$ the function $\tilde f(x)$ is continuous (moreover, it is two-times differentiable).
It follows from (\ref{leb2}) that  $P(f(W_t)=\tilde f(W_t))=1$ for any $t\ge0$ and since $Ef(W_t)=E\tilde f(W_t)$,
we obtain that
the processes $M_t={f(W_t)}-{Ef(W_t)}$ and $\tilde M_t={\tilde f(W_t)}-{E\tilde f(W_t)}$ are equivalent, which implies that
 the process $\tilde M_t$
is a continuous martingale. Therefore, it follows from part a) of this theorem that $\tilde f(x)$ is of the form (\ref{form23}) and hence,
$f(x)$ coincides with  the function (\ref{form23}) almost everywhere with respect to the Lebesgue measure.

The converse is proved similarly to the part a) of this theorem.
\end{proof}

{\bf Corollary 1.} Let $f=(f(x), x\in R)$ be a function of one variable.

a)  If  the process $\left(f\left(W_{t}\right), {\cal F}_t, t\geq 0\right)$  is  a  right-continuous
martingale,  then
\begin{equation}\label{abx}
f\left(x\right)=bx+c\;\;\;\text{for all}\;\;\;x\in R
\end{equation}
for some constants $b, c$.

b) If the process $\left(f\left(W_{t}\right), {\cal F}_t, t\geq 0\right)$  is a martingale, then $f\left(x\right)=bx+c$
almost everywhere with respect to the Lebesgue measure for some constants $b, c\in R$.
\begin{proof}
 If the process $f(W_t)$ is a martingale, then $g(t)=Ef(W_t)$ is  constant and the coefficient $a$ in (\ref{f4}) is equal to zero.
Therefore,  this corollary follows from Theorem 1.
\end{proof}

\

\section{Application to  Functional Equations.}

It was proved in \cite{M} (see also \cite{MT2} for multidimensional case)  that  if the function $f=(f(x),x\in R)$ is a measurable solution of the Cauchy additive functional equation
$$
f(x+y)=f(x)+f(y),\;\;\;\text{for all}\;\;\;x,y \in R,
$$
then the transformed process $(f(W_t), t\ge0)$ is a right-continuous martingale, which (by Corollary 1 of Theorem 2) implies that $f$ is a linear function.
Here we  propose a similar characterization of solutions of the quadratic functional equation
\begin{equation}\label{q1}
f\left(x+y\right)+f\left(x-y\right)=2f\left(x\right)+2f\left(y\right)\;\;\;\text{for all}\;\;\;x,y \in R,
\end{equation}
where for a measurable solution of this equation the process $f(W_t)$ is no longer a martingale, but  a martingale will be the difference $f(W_t) - Ef(W_t)$, which enable us to use Theorem 2
and to give a martingale characterization of equation (\ref{q1}).

It is well known (see, e.g.,  \cite{A2},\cite{SK}), that the general measurable solution of equation $(\ref{q1})$ is the function $f\left(x\right)=a x^{2}, a\in R$. Moreover, in
 \cite{A2} equation (\ref{q1}) has been solved without any assumptions using the Hamel basis. We consider only the Lebesgue measurable solutions of (\ref{q1})
and using Theorem 2 we characterize the general measurable solution of  equation (\ref{q1}) in terms of martingales, which gives also a probabilistic proof of this assertion.

Let first prove the following two lemmas.

{\bf Lemma 1.} Let  $f=(f(x), x\in R)$ be a measurable function and let the random variables $f(W_t)$ and $f(2W_s-W_t)$ be integrable. Then
for $s\le t$
\begin{equation}\label{lem}
E(f(2W_s-W_t)|F_s)=E(f(W_t)|F_s)\;\;\;\;P-\text{a.s.}
\end{equation}
In particular,
\begin{equation}\label{lem2}
Ef (W_t-2W_s)=Ef (2W_s-W_t)=Ef(W_t)\;\;\;\text{for }\;\;\;s\le t.
\end{equation}
\begin{proof}
Since  $F_t$ is generated by Brownian Motion $W_t$, it is sufficient to show that for any bounded measurable function $h(x_1, x_2,...,x_n)$ and
any set of indexes  $s_1, s_2,...,s_n\le s$
\begin{equation}\label{lem3}
E f(2W_s-W_t)h(W_{s_1},W_{s_2},..., W_{s_n})=E f(W_t)h(W_{s_1},W_{s_2},..., W_{s_n}).
\end{equation}
Since $E(2W_s-W_t)^2=E W_t^2=t$ and   $E(2W_s-W_t)W_{s_i}=E W_t W_{s_i}=s_i$, the distributions of normal vectors
$$
(W_t, W_{s_1},W_{s_2},..., W_{s_n})\;\;\text{and}\;\;\;(2W_s-W_t, W_{s_1},W_{s_2},..., W_{s_n})
$$
coincide, which implies equality
(\ref{lem3}).
\end{proof}

It is easy to show that if $f$ is a solution of (\ref{q1}), then  (see, e.g. \cite{SK})
\begin{equation}\label{r}
f(rx)=r^2f(x)\;\;\;\text{and}\;\;\;\; |f(rx)|=r^2|f(x)|
\end{equation}
for each rational $r$.

In particular, (\ref{r}) implies that
$$
f(r)=a r^2,\;\;\;\text{with}\;\;\;a=f(1)
$$
for any rational $r$.

It is evident that if $f$ is continuous then (\ref{r}) is satisfied for all real $r$, but we don't assume  continuity of the solution beforehand.

{\bf Lemma 2.} If $(f(x), x\in R)$ is a measurable solution of (\ref{q1}) then the random variable $f(\eta)$ is integrable for any random variable $\eta$ with normal distribution.
In particular, for a Brownian Motion $W$
\begin{equation}\label{w}
E|f(x+W_t)|<\infty, \;\;\;\text{for each}\;\;\; t\ge0\;\;\text{and}\;\;x\in R.
\end{equation}

{\it Proof.}  Let $r_n$ be sequence of rational numbers with $r_n\downarrow t,\;t\neq0$ and let $\xi$ be a random variable with standard normal distribution $(\xi\in N(0,1))$.
 By the  dominated convergence theorem we have
 $$
 E\f(\xi)e^{-r_n^2|f(\xi)|}\to E\f(\xi)e^{-t^2|f(\xi)|}\;\;\;\text{and}
$$
$$
E\f(\xi)e^{-|f(r_n\xi)|}=\int_R\f(\frac{1}{r_n}x)e^{-|f(x)|}\frac{e^{-\frac{x^2}{2r_n^2}}}{\sqrt{2\pi}r_n}dx
$$
$$
\to \int_R\f(\frac1{t}x)e^{-|f(x)|}\frac{e^{-\frac{x^2}{2t^2}}}{\sqrt{2\pi}t}dx=E\f(\xi)e^{-|f(t\xi)|},
$$
for each bounded continuous function $\f$.
Therefore, by (\ref{r}) these limits should coincide
 $E\f(\xi)(e^{-t^2|f(\xi)|}-e^{-|f(t\xi)|})=0$  and by arbitrariness of $\f$ we get
\begin{equation}\label{t}
|f(t\xi)|=t^2|f(\xi)|,\;\;P-a.s.\;\;\;\text{ for each}\;\;\;t\in R.
\end{equation}
Let $\eta$ be a gaussian random variable with mean $\mu$ and variance $\sigma^2$, independent from $\xi$. Then, it follows from (\ref{t}) that
\[
E\left(e^{-\eta^2|f(\xi)|}-e^{-|f(\eta\xi)|}\right)^2=
\int_R\frac{e^{-\frac{(y-\mu)^2}{2\sigma^2}}}
{\sqrt{2\pi}\sigma}
E\left(e^{-y^2|f(\xi)|}-e^{-|f(y\xi)|}\right)^2dy=0.
\]
Thus $\eta^2{|f(\xi)|}={|f(\eta\xi)|}$ and $(\mu^2+\sigma^2)|f(\xi)|=E(|f(\eta\xi)|/\xi)$ $P$-a.s..
Finally we get
\begin{equation}\label{fin}
E|f(x\eta)|=(\mu^2+\sigma^2)|f(x)|<\infty,\;\;\;\text{-a.e.}
\end{equation}
 with respect to the Lebesgue measure. This also implies that
$E|f(\eta)|<\infty$.

Indeed,  it follows from (\ref{fin}) that there exists $\gamma>1$ such that
\begin{equation}\label{g}
E|f(\gamma\eta)|<\infty.
\end{equation}
Therefore, after changing densities and taking the maximum
in the exponent we obtain from (\ref{g}) that
$$
E|f(\eta)|=E|f(\gamma\eta)|\gamma e^{-\frac{(\gamma\eta-\mu)^2}{2\sigma^2} +\frac{(\eta-\mu)^2}{2\sigma^2}}\le
$$
$$
\le\gamma e^{\frac{\mu^2}{2\sigma^2}\frac{\gamma-1}{\gamma+1}}E|f(\gamma \eta)|<\infty.
$$

{\bf Remark 1. } Lemma 2 implies that  any measurable solution of (\ref{q1}) is locally integrable.
Similar assertion for Cauchy's additive functional equation was proved in \cite{SM} using the Bernstein theorem on characterization of the normal distributions.

{\bf Theorem 3}.
The following assertions are equivalent:

i) the function  $f=(f(x), x\in R)$  is a measurable solution of (\ref{q1}),

ii)  $f=(f(x), x\in R)$ is a measurable even function with $f(0)=0$  and such that
$f(W_t)$ is integrable for every $t$ and the  process
$$
N_t={f(W_t)} - {Ef(W_t)},\;\;\; t\ge0,
$$
is a  right-continuous martingale,

iii) the function $f$ is of the form
\begin{equation}\label{q2}
f(x)= a x^2 ,
\end{equation}
for some constant $a\in R$.

\begin{proof}
$i)\to ii)$  It is evident that if $f$ is a solution of  $(\ref{q1})$ then $f\left(0\right)=0$ and
$f\left(x\right)=f\left(-x\right)$ for all $ x\in R$.    Therefore,
$$
f(W_t-2W_s-x)=f(2W_s-W_t+x)
$$
and substituting $x=W_t-W_s$ and $y=x+W_s$ in equation (\ref{q1}) we have that
\begin{equation}\label{q3}
f(x+W_t)+f(2W_s-W_t+x)=2f(W_t-W_s)+ 2f(x+W_s).
\end{equation}
Since $E(W_t-2W_s)^2=t$ for $s\le t$, the random variables $W_t$ and $2W_s-W_t$ have the same normal distributions and
by Lemma 2 we may take expectations in (\ref{q3}) to obtain
\begin{equation}\label{q5}
Ef(x+W_t)=Ef(W_t-W_s)+Ef(x+W_s)\;\;\;\;s\le t.
\end{equation}

If we take $s=0$ in (\ref{q5}) we get
$$
f(x)=E f(x+W_t)-Ef(W_t)=\int_R f(x+y) \frac{1}{\sqrt{2\pi t}}e^{-\frac{y^2}{2t}}dy-Ef(W_t)=
$$
 \begin{equation}\label{cont}
=\int_R f(y) \frac{1}{\sqrt{2\pi t}}e^{-\frac{(y-x)^2}{2t}}dy-Ef(W_t),
\end{equation}
which implies that $f$ is continuous by (\ref{w}).

Taking now conditional expectations in (\ref{q3}) for $x=0$, using the independent increment property of $W$ and then equality (\ref{q5})
we have that $P$-a.s.
$$
E(f(W_t)|F_s)+E(f(W_t-2W_s)|F_s)=2E(f(W_t-W_s)|F_s)+ 2f(W_s)=
$$
\begin{equation}\label{q6}
=2Ef(W_t-W_s)+ 2f(W_s)=2Ef(W_t)+2f(W_s)- 2Ef(W_s).
\end{equation}
On the other hand, it follows from Lemma 1 that $P$-a.s.
\begin{equation}\label{q7}
E(f(W_t-2W_s)|F_s)=E(f(W_t)|F_s).
\end{equation}
Therefore the martingale equality
$$
E(f(W_t)-Ef(W_t)|F_s)=f(W_s)-EfW_s),\;\;\;P-\text{a.s.},
$$
follows from equations (\ref{q6}) and (\ref{q7}). Thus, the process $N=({f(W_t)} - {Ef(W_t)}, t\ge0)$ is a martingale
with $P$-a.s. continuous  paths.

$ii)\to iii)$.  Since ${f(W_t)} - {Ef(W_t)}$ is a martingale, Theorem 2 implies that  the function $f$ should be of the form
$$
f(x)= ax^2+bx+c,\;\;\;a, b, c \in R.
$$
Since  $f$ is even we have that  $b=0$ and $c=0$  since $f(0)=0$. Thus, $f(x)=ax^2$ for some $a\in R$.

The proof of implication $iii)\to i)$ is evident.
\end{proof}

\

Now we give an application of Theorem 1.

Let consider the D'Alembert functional equation
\begin{equation}\label{dq1}
f\left(x+y\right)+f\left(x-y\right)=2f\left(x\right) f\left(y\right)\;\;\;\text{for all}\;\;\;x,y \in R,
\end{equation}
This equation possesses the following  solutions and only these:
$f(x)=0,   f(x) = \cosh\lambda x, f(x) = \cos\lambda x,$ where $\lambda$ is some constant.
The last two also contain (for $\lambda=0)$
the constant solution $f(x)=1$ (see, e.g.,  \cite{A}, \cite{SK}).

In the following theorem we give a martingale characterization of measurable strictly positive solutions of equation (\ref{dq1}).

{\bf Theorem 4}. The following assertions are equivalent:

i) the function  $f=(f(x), x\in R)$  is a measurable strictly positive solution of (\ref{dq1}),

ii)  $f=(f(x), x\in R)$ is a strictly positive even function with $f(0)=1$  and such that
$f(W_t)$ is integrable for every $t$ and the  process
$$
N_t=\frac{f(W_t)}{Ef(W_t)},\;\;\; t\ge0,
$$
is a  right-continuous martingale,

iii) the function $f$ is of the form
\begin{equation}\label{form2}
f(x)= \cosh(\lambda x) =\frac{1}{2}(e^{\lambda x}+e^{-\lambda x}) ,
\end{equation}
for some constant $\lambda\in R$.
\begin{proof} $i)\to ii).$
With $y=0$, it follows from (\ref{dq1}) that $f(x)=f(x)f(0)$, which implies $f(0)=1$, since we consider only solutions with $f(x)>0$. It is also evident that
$f$ is an even function, since taking $x=0$ from (\ref{dq1}) we have $f(y)+f(-y)=2f(0)f(y)$, hence $f(y)=f(-y)$.

Substituting $x=W_t-W_s$ and $y=x+W_s$ in equation (\ref{dq1}) we have that
\begin{equation}\label{dq3}
f(x+W_t)+f(W_t-2W_s-x)=2f(W_t-W_s)f(x+W_s).
\end{equation}
Since  $f(x)$ is positive, expectations bellow have a sense and  using the independent increment property of $W$
and  Lemma 1, we obtain from (\ref{dq3}) that
\begin{equation}\label{dq5}
Ef(x+W_t)=Ef(W_t-W_s)Ef(x+W_s)\;\;\;\;s\le t.
\end{equation}
Let $g(t)=Ef(W_t)$. Then $Ef(W_t-W_s)=g(t-s)$ and it follows from (\ref{dq5})  that $g$ satisfies the Cauchy exponential functional equation
$$
g(t)=g(t-s)g(s),\;\;\;s\le t
$$
on $R^+$. As it is well known (see, e.g., \cite{A}) that any bounded from bellow  solution of this equation
is of the form
\begin{equation}\label{aq7}
g(t)=e^{c t}\;\;\;\;\text{ for some
constant}\;\;\;c\in R.
\end{equation}
Therefore, $f(W_t)$ is integrable for any $t\ge0$ and $Ef(W_t)=e^{ct}$.

If we take $s=0$ in (\ref{dq5}) we obtain
$$
f(x)=\frac{E f(x+W_t)}{Ef(W_t)}=\frac{1}{Ef(W_t)}\int_R f(x+y) \frac{1}{\sqrt{2\pi t}}e^{-\frac{y^2}{2t}}dy=
$$
 \begin{equation}\label{cont}
=\frac{1}{Ef(W_t)}\int_R f(y) \frac{1}{\sqrt{2\pi t}}e^{-\frac{(y-x)^2}{2t}}dy,
\end{equation}
which implies that $f$ is continuous, since  $f(W_t)$ is integrable.

Taking now conditional expectations in (\ref{dq3}) for $x=0$, using the independent increment property of $W$ and then equality (\ref{dq5})
we have that $P$-a.s.
$$
E(f(W_t)|F_s)+E(f(W_t-2W_s)|F_s)=2E(f(W_t-W_s)|F_s)f(W_s)=
$$
\begin{equation}\label{dq6}
=2f(W_s) Ef(W_t-W_s)=2f(W_s)\frac{Ef(W_t)}{f(W_s)}.
\end{equation}
On the other hand, it follows from Lemma 1 that $P$-a.s.
\begin{equation}\label{dq7}
E(f(W_t-2W_s)|F_s)=E(f(2W_s-W_t)|F_s)=E(f(W_t)|F_s).
\end{equation}
Therefore the martingale equality
$$
E(\frac{f(W_t)}{Ef(W_t)}|F_s)=\frac{f(W_s)}{EfW_s)},\;\;\;P-\text{a.s.},
$$
follows from equations (\ref{dq6}) and (\ref{dq7}). Thus, the process $N=(\frac{f(W_t)}{Ef(W_t)}, t\ge0)$ is a martingale
with $P$-a.s. continuous  paths.

$ii)\to iii)$ It follows from Theorem 1 that $f(x)$ is of the form
$$
f(x)=a e^{\lambda x}+b e^{-\lambda x}
$$
for some constants $a, b, \lambda\in R$. Since the function $f$ is even we have that $a=b$ and $a+b=1$ by equality  $f(0)=1$. Therefore, $a=b=1/2$,
which implies representation (\ref{form}).

The proof of implication $iii)\to i)$ is evident.

\end{proof}

\section{On Time-dependent Martingale Transformations of a Brownian Motion.}

In this section we consider time-dependent functions $(f(t,x), t\ge0,x\in R)$ for which the transformed  processes
\begin{equation}\label{tr}
{f(t, \sigma W_t)}-{Ef(t,\sigma W_t)}\;\;\;\text{ and}\;\;\;{f(t, \sigma W_t)}/{Ef(t,\sigma W_t)}
\end{equation}
are martingales, where $\sigma$ is a constant. To obtain simple structural properties for such functions, as for the case of functions $f=(f(x), x\in R)$
in section 2, we need some type
of growth conditions on the function $f$, or one should  require the martingale property for transformed processes (\ref{tr}) for at least two different $\sigma\neq 0$.

Recall that a heat polynomial is any polynomial solution
of the heat equation $u_t+\frac12u_{xx}=0$.

{\bf Theorem 5.} Let $f:(0,\infty)\times R\to R $ be a continuous function, such that  $f(t,W_t)$ is integrable for every $t\ge0$ and $n$ is an integer $n\ge 1$. Then
 $f(t,W_t)-Ef(t,W_t)$ is a martingale
satisfying condition:

$C)$ for some $C>0$ the process
$$
 \la f(\cdot,W)-g\ra_t-C \int_0^t(1+s+W_s^2)^{n-1}ds\;\;\text{is  non-increasing},
$$
if and only if
the function $f(t,x)$ is of the form
\begin{equation}\lbl{qvdr}
f(t, x)=P(t,x)+c(t)
\end{equation}
for some heat polynomial $P$ of degree $n$ and deterministic function $c(t), t\ge0$.
\begin{proof}
Let us prove the sufficiency, the  necessity  part of this theorem  is evident.
The martingale property of the process  $f(t,W_t)-g(t)$ and the continuity of the function $f$ imply that
\begin{equation}\lbl{mart}
f(t, x)-g(t)=\int_R (f(T,y)-g(t))\frac{1}{\sqrt{2\pi(T-t)}}e^{-\frac{(y-x)^2}{2(T-t)}}dy.
\end{equation}

It is evident that $f-g$ is a  weak solution of the heat equation. Indeed for every infinitely differentiable finite (on $(0,\infty)\times {R}$) function $\varphi$, with
$supp\vp\subset(0,T)\times R$ for some $T>0$ we have from (\ref{mart}) that
\beaa
\int_0^\infty\int_{R} (f(s,y)-g(s))\big(\frac{\partial\varphi}{\partial t}(s,y)-\frac{1}{2}\varphi_{yy}(s,y)\big)dyds\\
=\int_0^T\int_{R} \int_R(f(T,x)-g(T))\frac{e^{-\frac{(x-y)^{2}}{2(T-s)}}}{\sqrt{2\pi (T-s)}}dx\big(\frac{\partial\varphi}{\partial t}(s,y)-\frac{1}{2}\varphi_{yy}(s,y)\big)dyds\\
=-\int_0^T\int_{R}(f(T,x)-g(T))\int_R\vp(s,y)\big(\frac{\partial}{\partial t}+\frac{1}{2}\frac{\partial^2}{\partial y^2}\big)\frac{e^{-\frac{(x-y)^{2}}{2(T-s)}}}{\sqrt{2\pi (T-s)}}dydxds=0.
\eeaa
By hypoelipticity property of the heat equation $f(t,x)-g(t)$ coincides with an infinitely differentiable function a.e. . Since $f(t,x)$ is continuous the function $f(t,x)-g(t)$ itself will be
infinitely differentiable. Therefore,
by the Ito formula we get
$$f(t,W_t)-g(t)= f(0,0)-g(0)+\int_0^tf_x(s,W_s)dW_s,$$
which implies that $\la f(\cdot,W)-g\ra_t=\int_0^t( f_x(s,W_s))^2ds$.
The condition of  this theorem gives  $|f_x(s,W_s)|^2\le C(1+s+W_s^2)^{n-1}$,  which is equivalent to the inequality $|f_x(t,x)|\le C(\sqrt{1+t+x^2})^{{n-1}}$.
Since $f_x$ is the classical  solution of the heat equation
by the Liouville theorem for heat equations \cite{HMa} it follows that $f_x(t,x)$ coincides a.e. with a heat polynomial of degree $n-1$. Therefore $f(t,x)$ is of the form (\ref{qvdr}).
\end{proof}

{\bf Corollary 2.} Let $f:(0,\infty)\times R\to R $ be a continuous function. Then
 $f(t,W_t)-Ef(t,W_t)$ is a martingale
satisfying condition:
$$
\;\text{the process}\;\; \la f(\cdot,W)-g\ra_t-C \int_0^t(1+s+W_s^2)ds\;\;\text{is  non-increasing for some}\;\;C>0,
$$
if and only if
the function $f(t,x)$ is of the form
\begin{equation}\lbl{qvdr2}
f(t, x)=a x^2+bx+c(t)
\end{equation}
for some constants $a, b\in R$ and deterministic function $c(t)=f(t,0), t\ge0$.

{\bf Remark 2.} If we only assume that the function $f(t,x)$ is measurable (without assumption of continuity of $f$), then the equations (\ref{qvdr}) and  (\ref{qvdr2}) will be satisfied almost
everywhere with respect to the Lebesgue measure $dt\times dx$.

{\bf Remark 3.} For each probability measure on $R$
the process $f(t,W_t)=\int_Re^{\sg W(t)-\frac{\sg^2}2t}d\nu(\sg)$ is a martingale and $g(t)=1$. Hence the martingale function is not even polynomial
without some type growth condition on $\la f(\cdot,W)-g\ra$ .

{\bf Remark 4.} One can prove that each heat polynomial of degree $n$ admits decomposition $f(t,x)=\sum_{k=0}^nC_kH_k(t,x)$ with respect to Hermite polynomials
$H_k(t,x)={(-t)^k}e^{x^2/2t}\frac{\partial^k}{\partial x^k}e^{-x^2/2t}$ (see the Appendix). Hence instead of ({\ref{qvdr}) we can write that there exist constants $C_k,k=0,...,n$ and  a function $c(t)$ such that
\beaa
f(t, x)=\sum_{k=0}^nC_kH_k(t,x)+c(t).
\eeaa

Now we give another description of the time-dependent martingale functions.

{\bf Theorem 6}. Let $f=(f(t, x), t\ge0, x\in R)$ be a  continuous   function. The following assertions are equivalent:

a) $E|f(t,\sigma W_t)|<\infty$ for every $t\ge0$ and  the  process
$$
M_t(\sigma)={f(t, \sigma W_t)}-{Ef(t,\sigma W_t)},\;\;\; t\ge0,
$$
is a  martingale for all $\sigma\in R$,

b)  $E|f(t,\sigma_1 W_t)|<\infty, E|f(t,\sigma_2 W_t)|<\infty$ for every $t\ge0$ and the processes
$$
M_t(\sigma_1)={f(t, \sigma_1 W_t)}-{Ef(t,\sigma_1 W_t)},\;\;\;  M_t(\sigma_2)={f(t, \sigma_2 W_t)}-{Ef(t,\sigma_2 W_t)}  \;\;\; t\ge0,
$$
 are martingales
 for two different $\sigma_1\neq\sigma_2\neq0.$

c) the function $f(t,x)$ is of the form
\begin{equation}\label{form25}
f(t, x)=a x^2+bx+c(t)
\end{equation}
for some constants $a, b\in R$ and deterministic function $c(t)=f(t,0), t\ge0$.
\begin{proof}
$b)\to c)$. Let
$$
g(t,\sigma)=Ef(t,\sigma W_t)\;\;\;\text{and}\;\;\;u(t,x)=f(t,x)-g(t,\sigma).
$$
Since $u(t,x)$  is continuous and the process $u(t,W_t)$ is a martingale,
 $u(t,x)$ will be of the class $C^{1.2}$ on $(0,T)\times R$ and  satisfies the  ''backward'' heat  equation
 (see, e.g. \cite{KAR})
\begin{equation}\label{heat}
\frac{\partial u}{\partial t}+\frac{\sigma^2}{2}\frac{\partial^2u}{\partial x^2}=0,\;\;\;\;0<t<T, x\in R,
\end{equation}
which implies that
\begin{equation}\label{heat2}
\frac{\partial\big(f(t, x)-g(t,\sigma)\big)}{\partial t}+\frac{\sigma^2}{2} \frac{\partial^2f(t,x)}{\partial x^2}=0,\;\;\;\;0<t<T, x\in R.
\end{equation}
Taking the difference of equations (\ref{heat2}) for $\sigma_1$ and $\sigma_2$ we have that
\begin{equation}\label{heat3}
\frac{\sigma_1^2-\sigma_2^2}{2}  \frac{\partial^2f(t,x)}{\partial x^2}=\frac{\partial\big(g(t,\sigma_1)-g(t,\sigma_2)\big)}{\partial t},\;\;\;\;0<t<T, x\in R.
\end{equation}
It follows from the last equation that the second derivative  $f_{xx}(t,x)$ is constant for any fixed $t$, which implies that $f(t,x)$ is a square trinomial with time dependent coefficients $a(t), b(t), c(t), t\ge0$
\begin{equation}\label{heat4}
f(t, x)=a(t) x^2+b(t)x+c(t).
\end{equation}
Therefore,
\begin{equation}\label{heat5}
g(t,\sigma)=Ef(t,\sigma W_t)=\sigma^2 t a(t)+c(t)
\end{equation}
and substituting expressions  (\ref{heat4}) and (\ref{heat5}) for $\sigma_1$ and $\sigma_2$ in (\ref{heat3}) we obtain that
$$
(\sigma_1^2-\sigma_2^2)a(t)=\frac{\partial \big((\sigma_1^2-\sigma_2^2)t a(t)\big)}{\partial t},
$$
which implies that $a'(t)t=0$ and hence $a(t)$ is a constant for any  $t>0$.

This, together with (\ref{heat4}) and (\ref{heat5}, implies that
$$
f(t,\sigma W_t)-Ef(t,\sigma W_t)=a\sigma^2(W_t^2-t)+\sigma b(t)W_t
$$
and this process is a martingale if and only if $b(t)$ is equal to a constant. Thus, $f(t,x)$ will be of the form (\ref{form25}.

$c)\to a)$ If  the function $f(t,x)$ is of the form (\ref{form25})
then
$$
f(t, \sigma W_t)= a\sigma^2 W_t^2+b\sigma W_t+c(t),\;\;\;\;
$$
 $E|f(t,\sigma W_t)|<\infty$ for all $t\ge0, \sigma\in R$ and  $Ef(t,\sigma W_t)=a\sigma^2 t+c(t)$.
It is evident that the process $f(t,\sigma W_t)-Ef(t,\sigma W_t)=a\sigma^2(W_t^2-t)$ is a martingale.

The implication $a)\to b)$ is evident.
\end{proof}

{\bf Corollary 3.} Let $f=(f(t, x), t\ge0, x\in R)$ be a  continuous   function. The following assertions are equivalent:

a) the  process  ${f(t, \sigma W_t)}\;\;\; t\ge0,$
is a  martingale for all $\sigma\neq 0$,

b)  the processes ${f(t, \sigma_1 W_t)}$ and ${f(t, \sigma_2 W_t)}$ are martingales  for two different
$\sigma_1\neq\sigma_2\neq0,$

c) the function $f(t,x)$ is of the form
\begin{equation}\label{form26}
f(t, x)=bx+c
\end{equation}
for some constants $b$ and  $c$.

\begin{proof}
 If the process $f(t, \sigma W_t)$ is a martingale, then $g(t)=Ef(t, \sigma W_t)$ is  constant and  from (\ref{heat5})
 $a(t)=a=0$ and $ c(t)$ is equal to a constant.
Therefore,  this corollary follows from Theorem 3.
\end{proof}

{\bf Theorem 7}. Let $f=(f(t, x), t\ge0, x\in R)$ be a  continuous  strictly positive  function differentiable at $t$ for any $t\ge0$. The following assertions are equivalent:

a)  $E|f(t,\sigma W_t)|<\infty$  for every $t\ge0$ and  the  process
$$
N_t(\sigma)=\frac{f(t, \sigma W_t)}{Ef(t,\sigma W_t)},\;\;\; t\ge0,
$$
is a  martingale for all $\sigma\in R.$

b)   $E|f(t,\sigma_1 W_t)|<\infty, E|f(t,\sigma_2 W_t)|<\infty$  for every $t\ge0$ and the processes
$$
N_t(\sigma_1)=\frac{f(t, \sigma_1 W_t)}{Ef(t,\sigma_1 W_t)},\;\;\; N_t(\sigma_2)=\frac{f(t, \sigma_2 W_t)}{Ef(t,\sigma_2 W_t)},\;\;t\ge0,
$$
are martingales  for two different
 $\sigma_1\neq\sigma_2\neq0.$

c) the function $f(t,x)$ is of the form
\begin{equation}\label{fform25}
f(t, x)=a c(t) e^{\lambda x}+b c(t) e^{-\lambda x}
\end{equation}
for some constants $a\ge0, b\ge0$ with $a+b=1, ab\neq0$ and deterministic function $c(t)=f(t,0), t\ge0$.
\begin{proof}
$b)\to c)$. Let
$$
g(t,\sigma)=Ef(t,\sigma W_t)\;\;\;\text{and}\;\;\;h(t,x)=\frac{f(t,x)}{g(t,\sigma)}.
$$
Since $h(t,x)$ is continuous and the process $h(t,W_t)$ is a martingale,
 $h(t,x)$ will be of the class $C^{1.2}$ on $(0,T)\times R$ and  satisfies the  ''backward'' heat  equation
 (see, e.g. \cite{KAR})
\begin{equation}\label{fheat}
\frac{\partial h}{\partial t}+\frac{\sigma^2}{2}\frac{\partial^2h}{\partial x^2}=0,\;\;\;\;0<t<T, x\in R,
\end{equation}
which implies that
\begin{equation}\label{fheat2}
\frac{\partial\big(\frac{f(t, x)}{g(t,\sigma)}\big)}{\partial t}+\frac{\sigma^2}{2} \frac{\partial^2f(t,x)}{\partial x^2}\frac{1}{g(t,\sigma)}=0,\;\;\;\;0<t<T, x\in R.
\end{equation}
Since $f(t,x)$ is differentiable at $t$, the function $g(t,\sigma)$ will be also differentiable and from (\ref{fheat2}) we have that
\begin{equation}\label{fheat3}
\frac{\sigma^2}{2} \frac{\partial^2f(t,x)}{\partial x^2}+\frac{\partial\big({f(t, x)}\big)}{\partial t}-f(t,x)\frac{g'(t,\sigma)}{g(t,\sigma)}=0,\;\;\;\;0<t<T, x\in R.
\end{equation}

Taking the difference of equations (\ref{fheat3}) for $\sigma_1$ and $\sigma_2$ we have that
\begin{equation}\label{fheat4}
\frac{\sigma_1^2-\sigma_2^2}{2}  \frac{f_{xx}(t,x)}{f(t,x)}=\frac{g'(t,\sigma_1)}{g(t,\sigma_1)}-\frac{g'(t,\sigma_2)}{g(t,\sigma_2)},\;\;\;\;0<t<T, x\in R.
\end{equation}
It follows from the last equation that  $f_{xx}(t,x)/f(t,x)$ does not depend on $x$, i.e.,  ${f_{xx}(t,x)}/{f(t,x)}=c(t)$ for some function $(c(t),t\ge0)$, which should be positive for all $t\ge0$, since
if $c(t_0)<0$ for some $t_0$ then the general solution of equation ${f_{xx}(t_0,x)}/{f(t_0,x)}=c(t_0)$ leads to $f(t_0,x)$ which can take negative values. Hence
\begin{equation}\label{fheat6}
  \frac{f_{xx}(t,x)}{f(t,x)}=\lambda^2(t),
\end{equation}
for some function $(\lambda(t),t\ge0)$. For any fixed $t$ the general solution of equation (\ref{fheat6}) is of the form
\begin{equation}\label{fheat7}
  f(t,x)=a(t)e^{\lambda(t)x}+b(t)e^{-\lambda(t)x},
\end{equation}
for some functions of $t$ - $a(t), b(t)$ and $\lambda(t)$.

Let first show that $\lambda(t)=\lambda$ for all $t\ge0$, for some $\lambda\in R$.

It follows from (\ref{fheat7})
\begin{equation}\label{fheat8}
 E f(t,\sigma W_t)=(a(t)+b(t))e^{\frac{\sigma^2\lambda^2(t)}{2}t}
\end{equation}
and it is easy to see that
\begin{equation}\label{fheat9}
\frac{g'(t,\sigma)}{g(t,\sigma)}=\frac{a'(t)+b'(t)}{a(t)+b(t)}+\sigma^2\lambda(t)\lambda'(t) t+\frac{\sigma^2}{2}\lambda^2(t).
\end{equation}
Substituting expressions  (\ref{fheat9})  for $\sigma_1$ and $\sigma_2$ in (\ref{fheat4}) we obtain from (\ref{fheat6}) that
$\lambda(t)\lambda'(t) t=0$, which implies that
\begin{equation}\label{fheat90}
\lambda^2(t)=\lambda^2\;\;\;\text{ for some constant}\;\;\;\lambda\in R.
\end{equation}
Therefore, it follows from (\ref{fheat90}),  (\ref{fheat7}) and  (\ref{fheat8}) that
\begin{equation}\label{fheat10}
N_t(\sigma)=\frac{f(t, \sigma W_t)}{Ef(t,\sigma W_t)}=\frac{a(t)}{a(t)+(b(t)}e^{\lambda\sigma W_t-\frac{\sigma^2\lambda^2}{2}t}+
\frac{b(t)}{a(t)+(b(t)}e^{-\lambda\sigma W_t-\frac{\sigma^2\lambda^2}{2}t}.
\end{equation}
Since the processes $X_t=e^{\lambda\sigma W_t-\frac{\sigma^2\lambda^2}{2}t}$ and  $Y_t=e^{-\lambda\sigma W_t-\frac{\sigma^2\lambda^2}{2}t}$ are martingales and
$P(X_t\neq Y_t)=1$ for all $t$,
the process $N_t(\sigma)$ defined by (\ref{fheat10}) will be a martingale if and only if
\begin{equation}\label{fheat11}
\alpha_t\equiv\frac{a(t)}{a(t)+(b(t)}=a,\;\;\;\text{and}\;\;\;\beta_t\equiv\frac{b(t)}{a(t)+(b(t)}=b
\end{equation}
for some constants $a, b\in R$.

Indeed, since $\alpha_t$  is a deterministic function, $\beta_t=1-\alpha_t$ and the process $\alpha_t X_t+(1-\alpha_t)Y_t$ is a martingale, then for any $s\le t$
 $$
 \alpha_s X_s+(1-\alpha_s)Y_s=E(\alpha_t X_t+(1-\alpha_t)Y_t/F_s)=\alpha_t X_s+(1-\alpha_t)Y_s,
 $$
 which implies that
 $(\alpha_t-\alpha_s)(X_s-Y_s)=0$. Therefore, $\alpha_t=\alpha_s$ and $\alpha_t$ is equal to a constant by arbitrariness of $s$ and $t$.

Therefore,  (\ref{fheat90}),  (\ref{fheat7}) and (\ref{fheat11}) imply that
\begin{equation}\label{fheat12}
  f(t,x)=\big(a(t)+b(t)\big) [a e^{\lambda x}+be^{-\lambda x}],
\end{equation}
where by (\ref{fheat11}) $a+b=1$ and $c(t)\equiv a(t)+b(t)=f(t,0)$ from (\ref{fheat12}). Besides, $a\ge0, b\ge0$ and $ab\neq 0$, since $f(t,x)$ is strictly positive.
 Hence $f(t,x)$ is of the form (\ref{fform25}).

$c)\to a)$ If  the function $f(t,x)$ is of the form (\ref{fform25})
then
\begin{equation}\label{fheat13}
  f(t,\sigma W_t)=f(t,0) [a e^{\lambda\sigma W_t}+be^{-\lambda\sigma W_t}]
\end{equation}
 and $E|f(t,\sigma W_t)|<\infty$ for all $t\ge0, \sigma\in R$. It is evident that  $Ef(t,\sigma W_t)=f(t,0)e^{\frac{\sigma^2\lambda^2}{2}t}$
and the process
$$
\frac{f(t, \sigma W_t)}{Ef(t,\sigma W_t)}=a e^{\lambda\sigma W_t-\frac{\sigma^2\lambda^2}{2}t}
+ b e^{-\lambda\sigma W_t-\frac{\sigma^2\lambda^2}{2}t}
$$
 is a martingale for any $\sigma$.

The implication $a)\to b)$ is evident.
\end{proof}

\appendix

\section{Appendix}

The Hermite polynomial is defined by
$H_k(t,x)=(-t)^ke^{x^2/2t}\frac{\partial^k}{\partial x^k}e^{-x^2/2t}$. Using the Taylor expansion for exponential function we get
\beaa\sum_{k=0}^\infty\sg^k\frac{H_k(t,x)}{k!}
=\sum_{k=0}^\infty\sg^k\frac{(-t)^ke^{x^2/2t}\frac{\partial^k}{\partial x^k}e^{-x^2/2t}}{k!}\\
=e^{x^2/2t}\sum_{k=0}^\infty\frac{(-t\sg)^k}{k!}\frac{\partial^k}{\partial x^k}e^{-x^2/2t}=e^{x^2/2t}e^{-(x-\sg t)^2/2t}=e^{\sg x-\frac{\sg^2}2t}.
\eeaa
From the expansion $e^{\sg x-\frac{\sg^2}2t}=\sum_{n=0}^\infty\frac{H_n(t,x)}{n!}\sg^n$ follows that $\{H_n\}$ are  heat polynomials.
If $f(t,x)=\sum_{k=0}^na_k(t)x^{n-k}$ is arbitrary heat polynomial, then equalizing coefficients in the heat equation we get
\beaa
a_0'(t)=0,\;a_{2j}'(t)=-\frac{(n-2j+2)(n-2j+1)}{2}a_{2j-2}(t), \;0<2j\le n\\
a_1'(t)=0,\;a_{2j+1}'(t)=-\frac{(n-2j+3)(n-2j+2)}{2}a_{2j-1}(t),\;0<2j+1\le n.
\eeaa
Hence solution of the system
\beaa
a_0(t)=C_0,\\
a_{2j}(t,C_0,C_2,...,C_{2j})\\=-\frac{(n-2j+2)(n-2j+1)}{2}\int a_{2j-2}(t,C_0,C_2,...,C_{2j-2})dt+C_{2j}, \;0<2j\le n\\
a_1(t)=C_1,\\
a_{2j+1}(t,C_1,C_3,...,C_{2j+1})\\=-\frac{(n-2j+3)(n-2j+2)}{2}\int a_{2j-1}(t,C_1,C_3,...,C_{2j-1})dt+C_{2j+1},\;1<2j+1\le n.
\eeaa
linearly depends on arbitrary constants $C_0,C_1,...,C_n$. Since $H_k,k=0,1,...,n$ are linear independent heat polynomials with degree less than $n+1$, there exist
$C_0',...,C_n'$, such  that $f=\sum_{k=0}^nC_k'H_k.$

{\bf Example.} For $n=3$ and $f(t,x)=a_0(t)x^3+a_1(t)x^2+a_2(t)x+a_3(t)$
we get
\beaa
a_0=C_0,\; \;a_2(t)=-3C_0t+C_2,\\
a_1=C_1,\;\;a_3(t)=-C_1t+C_3,\\
f(t,x)=C_0x^3+C_1x^2+(-3C_0t+C_2)x-C_1t+C_3\\
=C_0(x^3-3tx)+C_1(x^2-t)+C_2x+C_3\\
=C_0H_3(t,x)+C_1H_2(t,x)+C_2H_1(t,x)+C_3H_0(t,x).
\eeaa

{}

\end{document}